\input amstex
\documentstyle{amsppt}
\magnification=1200
\loadeusm
\loadbold
\PSAMSFonts

\NoBlackBoxes
\nologo
\TagsOnRight
\rightheadtext{Harmonic Sections, and Cheeger-Gromoll Metrics}
\pagewidth{16truecm}
\pageheight{24.3truecm}
\linespacing{1.2}
\redefine\!{\kern-.075em}
\redefine\qed{\hfill$\ssize\square$}
\redefine\div{\operatorname{div}}
\predefine\accute{\'}
\redefine\'{\kern.05em{}}
\def\<{\langle}
\def\>{\rangle}
\def\leqs{\leqslant}
\def\geqs{\geqslant}
\def\volume{\operatorname{vol}}
\def\Tr{\operatorname{Trace}}
\def\E{\Cal E}
\def\bE{{\fam\cmbsyfam E}}

\def\C{\eusm C}
\def\F{\eusm F}
\def\G{\eusm G}
\def\W{\eusm W}
\def\a{\alpha}
\def\b{\beta}
\def\s{\sigma}
\def\bsigma{\boldsymbol\sigma}
\def\v{\varphi}
\def\upchi{\raise.4ex\hbox{$\chi$}}
\def\nab#1{\nabla\kern-.2em\lower.8ex\hbox{$\ssize#1$}\'}
\def\nabv#1{\nabla^v\kern-.4em\lower.8ex\hbox{$\ssize#1$}\'}

\def\dotprod{\'\raise.25ex\hbox{$\sssize\bullet$}\'}
\def\dvs{d^{\'v}\kern-.15em\s}
\topmatter
\title
Harmonic Sections of Riemannian Vector Bundles, and Metrics of
Cheeger-Gromoll Type
\endtitle
\author 
M. Benyounes,
E. Loubeau,
C. M. Wood
\endauthor
\address
D\accute epartement de Math\accute ematiques,
Laboratoire CNRS UMR 6205,
Universit\accute e de Bretagne Occidentale,
6 Avenue Victor Le Gorgeu, 
CS 93837,
29238 Brest Cedex 3,
France.
\endaddress
\email 
Michele.Benyounes\@univ-brest.fr,
Eric.Loubeau\@univ-brest.fr 
\endemail
\address 
Department of Mathematics, University of York, Heslington, York
Y010 5DD, UK.
\endaddress
\email 
cmw4\@york.ac.uk 
\endemail
\thanks
This research originated during the EDGE conference held at the
Universidad de Granada in February 2004, and the authors would like to
express their gratitude to the Departamento de Geometr\acuteaccent\i a y
Topolog\acuteaccent\i a for hosting and organizing the event.
\endthanks
\dedicatory
Dedicated to Professors J Eells and J H Sampson.
\enddedicatory
\abstract
We study harmonic sections of a Riemannian vector bundle $\E\to M$ whose
total space is equipped with a 2-parameter family of metrics $h_{p,q}$ which
includes both the Sasaki and Cheeger-Gromoll metrics.  The restrictions of
the $h_{p,q}$ to the total space of any sphere subbundle $S\E(k)$ of $\E$
(where $k>0$ is the radius) are essentially the same for all $(p,q)$, and
it is shown that for every $k$ there exists a unique $p$ such that the
harmonic sections of $S\E(k)$ are harmonic sections of $\E$ with respect to
$h_{p,q}$ for all $q$.  In both the compact and non-compact cases
Bernstein regions of the $(p,q)$-plane are identified, where the only
harmonic sections of $\E$ with respect to $h_{p,q}$ are parallel.  Examples
are constructed of compact vector fields which are harmonic sections of
$\E=TM$ in the case where $M$ has non-zero Euler characteristic.
\endabstract
\keywords
$(p,q)$-harmonic section, Sasaki metric, Cheeger-Gromoll metric, (strictly)
$q$-Riemannian section, Kato inequality, Bernstein region, Hopf vector
field, conformal gradient field
\endkeywords
\subjclass
53C43 (53C07, 53C24, 58E15, 58E20, 58G30)
\endsubjclass
\endtopmatter
\document
\head
1. Introduction
\endhead
The aim of this paper is to introduce some new criteria for deciding which
smooth vector fields on a smooth, oriented, connected (but not necessarily
compact) Riemannian manifold $(M,g)$, or in general which smooth sections
$\s$ of a smooth oriented Riemannian vector bundle $(\E,\<\,,\>,\nabla)\to
M$, qualify as ``better than the rest''.  In so doing we overcome some of
the limitations of existing criteria, which we briefly review.
\smallskip\noindent
{\bf1.}\quad
$\boldsymbol\nabla\bsigma\boldkey=\bold0$. Since the fibre metric $\<\,,\>$
is holonomy-invariant, and $M$ is connected, parallel sections have
constant length.  Therefore if the Euler class $\upchi(\E)\neq0$ there are
no non-trivial solutions.  (The trivial solution is of course the zero
section.)  Since the existence of solutions is equivalent to reduction of
the holonomy of $\nabla$, amongst the many other necessary conditions is de
Rham's decomposition theorem: if $\E=TM$ and $\<\,,\>=g$ then the universal
cover of $M$ splits as a Riemannian product $M'\times\Bbb R$.  So whilst
this criterion undeniably produces the ``best'' sections, its applicability
is severely limited.
\smallskip\noindent
(2)\quad
$\boldsymbol\Delta\bsigma\boldkey=\bold0$.  Here $\E=TM$, $\<\,,\>=g$,
$\nabla$ is the Levi-Civita connection, and $\Delta$ is the Hodge-de Rham
laplacian on $1$-forms, dualized to act on vector fields.  By Hodge's
theorem, if $M$ is compact then the solution space is isomorphic to
$H^1(M,\Bbb R)$; in particular if the first Betti number $\b_1(M)=0$ then
there are no non-trivial solutions.  Furthermore Bochner's vanishing
theorem informs us that when $M$ has positive Ricci curvature there are no
non-parallel solutions.
\smallskip\noindent
(3)\quad
$\bsigma$ {\bf is a harmonic section of} $\bE$ \cite{13,\,14}.  Here one
measures the  {\sl vertical energy\/} (or {\sl total bending\/} \cite{20})
of
$\s$:
$$
E^v(\s)\,=\,\frac12\int_M|\'\nabla\s\'|^2\,\volume(g),
\tag1-1
$$
(assuming for convenience that $M$ is compact; otherwise one works over
relatively compact domains), and looks for critical points with respect to
smooth variations through sections of $\E$.  The Euler-Lagrange equations
are once again linear:
$$
\nabla^*\nabla\s\,=\,0,
\tag1-2
$$
where $\nabla^*\nabla$ is the {\sl rough Laplacian:} 
$$
\nabla^*\nabla\,=\,-\Tr\nabla^2
$$
If $M$ is compact and (1-2) holds then integrating by parts:
$$
0\,=\,\int_M\<\nabla^*\nabla\s,\s\>\volume(g)\,
=\,\int_M|\'\nabla\s\'|^2\volume(g),
$$
so all harmonic sections of $\E$ are parallel.  The same is true if $M$ is
non-compact, provided $\s$ has constant length (see Lemma 3.4).
\smallskip\noindent
(4)\quad
$\boldkey|\'\bsigma\'\boldkey|\boldkey=\bold k$ {\bf (constant) and
$\bsigma$ is a harmonic section of the radius-$\bold k$ sphere bundle}
\cite{22,\,24}.  
Here the vertical energy functional (1-1) is restricted to sections of the
subbundle $S\E(k)\to M$, where:
$$
S\E(k)\,=\,\{e\in\E:|\'e\'|=k\},
$$
and the imposition of this constraint causes the Euler-Lagrange equations to
become mildly non-linear (see Remark 3.8):
$$
\nabla^*\nabla\s\,
=\,\frac{1}{k^2}\,|\'\nabla\s\'|^2\'\s
\tag1-3
$$
The solutions of (1-3) clearly include all parallel sections of length $k$
(if any), but when $\E=TM$ many additional solutions have been identified
\cite{1,\,8,\,9,\,17,\,19}, which in turn may be examined for stability
\cite{2,\,3,\,4,\,12,\,24}.  Unfortunately the theory is limited to bundles
with $\upchi(\E)=0$.
\medskip
Our new criteria remove the topological restriction $\upchi(\E)=0$, whilst
retaining all solutions of the constrained variational problem (4)
(Theorem A/4.1).  The basic idea is to obtain interesting non-linear
equations, such as (1-3), by altering the background metric data, rather
than introducing constraints.  We note first that definition (1-1) is
equivalent to:
$$
E^v(\s)\,=\,\frac12\int_M|\'\dvs\'|^2\volume(g),
\tag1-4
$$
where $\dvs$ is the vertical component of the differential $d\s$ with
respect to the connection $\nabla$ in $\E\to M$, and the norm in $T\E$ is
that of the {\sl Sasaki metric\/} $h$ on $\E$ \cite{18}.  The idea is to
study the functional (1-4) when $h$ is generalized to a $2$-parameter
family of metrics $h_{p,q}$ on $\E$, for which $h_{0,0}=h$ and $h_{1,1}$
is the {\sl Cheeger-Gromoll metric\/} \cite{7,\,15}.  (Both the Sasaki and
Cheeger-Gromoll metrics were originally defined for $\E=TM$, but generalize
in a natural way; see Remark 2.1.)  Other geometrically interesting metrics
occur in this family; for example $h_{2,0}$ is the stereographic metric
(Remark 2.2).
\par
Actually the term ``metric'' is used somewhat informally.  If $q\geqs0$ then
$h_{p,q}$ is indeed a Riemannian metric.  However if $q<0$ then $h_{p,q}$
has varying signature and is consequently not even semi-Riemannian: it is
Riemannian within the tubular neighbourhood of the zero section of radius
$1/\sqrt{-q}$, Lorentzian on the interior of the complement, and positive
semi-definite on the boundary.  This behaviour may be viewed as a
manifestation of {\sl Kato's inequality\/} \cite{6}.  A section whose image
lies in the closure of this tubular neighbourhood is said to be  {\sl
$q$-Riemannian\/} (see Remark 2.4).  If $q<0$ and $\s$ is not
$q$-Riemannian then it is possible that $E^v(\s)<0$.  In any case, if $\s$
is stationary for (1-4) with respect to the metric $h_{p,q}$ on $\E$, and
smooth variations through sections of $\E$ then we say that $\s$ is a  {\sl
$(p,q)$-harmonic section of $\E$.}  The Euler-Lagrange equations for
$(p,q)$-harmonic sections are derived in \S3 (Theorem 3.6), after a
somewhat lengthy sequence of calculations.  They are considerably more
complicated than (1-2), to which of course they reduce when $(p,q)=(0,0)$. 
However for all $(p,q)$ the parallel sections of $\E$ are always
$(p,q)$-harmonic, and amongst $q$-Riemannian sections they comprise the
absolute minima of $E^v$.
\par
An interesting feature of the $h_{p,q}$ is that they restrict to
essentially the same Riemannian metric on $S\E(k)$, even when $q<0$ and
$k>1/\sqrt{-q}\,$ (Remark 2.3).  Hence $(p,q)$-harmonic sections of $S\E(k)$
are characterised by equations (1-3) for all $(p,q)$, and may therefore be
referred to simply as {\sl harmonic sections of $S\E(k)$.}  For bundles
with $\upchi(\E)=0$ we establish the following relationship between
$(p,q)$-harmonic sections of $\E$ and harmonic sections of $S\E(k)$:
\proclaim{Theorem A}
Suppose that $|\'\s(x)|=k>0$ for all $x\in M$.
\flushpar
\rom{(a)}\quad 
If $p\neq1+1/k^2$ then $\s$ is a $(p,q)$-harmonic section of $\E$ if and
only if $\s$ is parallel.
\flushpar
\rom{(b)}\quad
If $p=1+1/k^2$ then $\s$ is a $(p,q)$-harmonic section of $\E$ if and only
if $\s$ is a harmonic section of $S\E(k)$.
\endproclaim
Theorem A may also be regarded as a first source of examples of
$(p,q)$-harmonic sections of $\E$, when $p>1$ and $\upchi(\E)=0$
(Example 4.2).  In seeking non-trivial examples of $(p,q)$-harmonic
sections of bundles with non-zero Euler class, we establish the following
rather more complicated set of restrictions on $(p,q)$:
\proclaim{Theorem B}
Suppose $M$ is compact, $\upchi(\E)\neq0$, and $\s$ is a non-trivial
section of $\E$.  For each $p\in\Bbb R$ there exists at most one 
$q\in\Bbb R$ such that $\s$ is a $(p,q)$-harmonic section of $\E$, and:
\flushpar
{\rm(a)}\quad
if $-4\leqs p\leqs-1$ then $\,q<-1-p$\rom;
\flushpar
{\rm(b)}\quad
if $-1\leqs p\leqs1$ then $\,q<0$\rom;
\flushpar
{\rm(c)}\quad
if $\,1<p\leqs2$ and $\,\|\'\s\'\|_\infty\leqs1/\sqrt{p-1}$ then
$\,q<0$\rom;
\flushpar
{\rm(d)}\quad
if $\,2\leqs p$ and $\,\|\'\s\'\|_\infty\leqs1/\sqrt{p-1}$ then
$\,q<1-p/2$.
\endproclaim
The appearance of $\|\'\s\'\|_\infty$ in Theorem B\,(c),(d) at first sight
seems counter-intuitive, since it implies that $(p,q)$-harmonicity is not
invariant under scaling (when $p>1$).  This reflects the non-linearity of
the $(p,q)$-harmonic section equations.  We do not know whether any
restrictions on the location of $q$ exist when $p<-4$.
\par
Theorem B is deduced from a Bernstein-type theorem (Theorem 4.6) and a
uniqueness theorem (Theorem 4.8), which in fact yield a more general result
valid for $(p,q)$-harmonic sections of non-constant length (Corollary
4.9).  Analogous results are available for the non-compact case (Theorem
4.3, Theorem 4.5, Corollary 4.10), although in order to compensate for the
unavailability of the global techniques used for Theorem B these are all
qualified by the assumption that $|\'\s\'|^2\colon M\to\Bbb R$ is a
harmonic function; they may therefore be viewed as generalizations of
Theorem A.  In both the compact and non-compact cases it transpires that no
additional examples of harmonic sections of $\E$ arise when the Sasaki
metric is replaced by the Cheeger-Gromoll metric \cite{16}, or indeed any
metric $h_{p,p}$ with $0\leqs p\leqs1$ (Remarks 4.4 and 4.7).  In fact all
metrics $h_{p.q}$ with $0\leqs p\leqs 1$ and $q\geqs0$  exhibit this
behaviour.  In order to find non-trivial examples of $(p,q)$-harmonic
sections of bundles with non-zero Euler class it is therefore necessary to
explore more ``remote'' regions of the $(p,q)$-plane.  
\par
The non-applicability of standard existence theory for harmonic maps (for
example \cite{10}), and its generalization to harmonic sections (for example
\cite{21,\,23}) necessitates a somewhat {\it ad hoc\/} approach to
the construction of examples.  In \cite{3} it was shown that normalizing a
{\sl conformal gradient field\/} on $S^5,S^7,\dots$ away from its (two)
zeroes produces a singular unit vector field whose energy infimizes the
energy functional when restricted to the space of smooth unit vector
fields.  (On $S^3$ the Hopf vector field is an absolute energy minimizer
\cite{4,\,12}.)  We show (Theorem 5.2) that if $\s$ is a conformal
gradient field on $M=S^n$ with $n\geqs3$ then $\s$ is a $(p,q)$-harmonic
section of $TM$ precisely when $p=n+1$, $q=2-n$, and
$\|\'\s\'\|_\infty=1/\sqrt{n-2}$.  (Note that these values of $(p,q)$ are
consistent with Theorem B.)  Although $q<0$, these $(p,q)$-harmonic sections
are $q$-Riemannian, but only just (Remark 5.3).  This example suggests that
in general a section should not be expected to be $(p,q)$-harmonic for more
than one metric $h_{p,q}$, although the existence of a $1$-parameter family
of $(p,q)$ is not precluded by Theorem B.  It also illustrates once again,
in a dramatic way, the non-invariance of solutions of the $(p,q)$-harmonic
section equations under scaling.  This suggests the following simple
general {\it ansatz:} given a ``trial'' section $\xi$, try to construct a
$(p,q)$-harmonic section by linearly rescaling $\xi$.  If this fails, try
a conformal rescaling.  Of course, the choice of $\xi$ remains {\it ad hoc.}
In this vein, we conclude by showing (Theorem 5.4) that when $M$ is an
odd-dimensional sphere the only $(p,q)$-harmonic sections of $TM$ obtained
by conformally rescaling the Hopf vector field $\xi$ are precisely those
covered by Theorem A: namely $\s=k\'\xi$ where $k=\pm1/\sqrt{p-1}\,$ and
$p>1$.
\bigskip
\head
2. The Vertical $(p,q)$-Energy Functional
\endhead
Let $(M,g)$ be a connected Riemannian $n$-manifold, and let 
$\pi\colon\E\to M$ be a vector bundle with connection $\nabla$ and
holonomy-invariant fibre metric $\<\,,\>$.
\subhead
Remark
\endsubhead
For the most part, connectedness of $M$ is simply a convenience which allows
us to simplify the exposition (for example, a parallel section of $\E$ then
has constant length), and most of our results are true whether or not $M$ is
connected.  The exceptions are Theorems 4.5\,(b) and 4.6\,(b) where
connectedness is an essential hypothesis.
\subhead
Remark
\endsubhead
The holonomy-invariance of $\<\,,\>$ will be used in many of our
calculations,  usually without comment, and is essential to our results.
\medskip
Let $K\colon T\E\to\E$ be the {\sl connection map\/} \cite{11} for $\nabla$:
$$
\CD
\E@<K<<T\E\\
@V\pi VV  @VV d\pi V \\
M@<<<TM
\endCD
$$
\medskip\noindent
and let $e\in\E$ and $A,B\in T_e\E$.  For any pair of parameters 
$p,q\in\Bbb R$ we define a symmetric $2$-covariant tensor $h_{p,q}$ on $\E$
as follows:
$$
h_{p,q}(A,B)\,=\,g(d\pi(A),d\pi(B))\,
+\,w^p(e)\bigl(\<KA,KB\>\,
+\,q\'\<KA,e\>\<KB,e\>\bigr),
\tag2-1
$$
where:
$$
w(e)\,=\,\frac{1}{1+|\'e\'|^2}
$$  
If $q\geqs0$ then $h_{p,q}$ is a Riemannian metric; however if $q<0$ then
$h_{p,q}$ is a Riemannian metric only on the following tubular
neighbourhood of the zero section:
$$
B\E(1/\sqrt{-q})\,
=\,\{e\in\E:|\'e\'|^2<-1/q\}
$$
\subhead
Remark 2.1
\endsubhead
If $(p,q)=(0,0)$ then $h_{p,q}$ is the {\sl Sasaki metric\/} \cite{18}:
$$
h_{0,0}(A,B)\,=\,g\bigl(d\pi(A),d\pi(B)\bigr)\,+\,\<KA,KB\>,
$$
whereas if $(p,q)=(1,1)$ then $h_{p,q}$ is the {\sl Cheeger-Gromoll
metric\/} \cite{7}:
$$
h_{1,1}(A,B)\,=\,g\bigl(d\pi(A),d\pi(B)\bigr)\,
+\,\frac{1}{1+|\'e\'|^2}\,\bigl(\<KA,KB\>
\,+\,\<KA,e\>\<KB,e\>\bigr)
$$
In all cases the bundle projection $(\E,h_{p,q})\to(M,g)$ is horizontally
isometric; in particular, if $q\geqs0$ it is a Riemannian submersion.
\subhead
Remark 2.2
\endsubhead
Another way of thinking about $h_{p,q}$ is as the horizontal lift of $g$ to
$\E$, supplemented by the metric on the fibres induced by the following
rotationally symmetric metric on Euclidean space:
$$
\frac{1}{(1+|\'x\'|^2)^p}\left(\sum_i(dx^i)^2\,
+\,q\sum_{i,j}x_i\,x_j\,dx^idx^j\right)
\tag2-2
$$
In particular, if $(p,q)=(2,0)$ then (2-2) is the stereographic metric,
up to homothety.
\subhead
Remark 2.3
\endsubhead
For any $k>0$, if $A,B$ are tangent to the total space of the sphere bundle
$S\E(k)\to M$ then $\<KA,e\>=0$ etc.  It therefore follows from equation
(2-1) that the restriction of $h_{p,q}$ to $S\E(k)$ is a Riemannian metric
for all $(p,q)$, and these metrics are essentially the same: the
restriction of $h_{p,q}$ differs from that of $h_{0,0}$ by the (constant)
factor $(1+k^2)^{-p}$ in the vertical direction.
\medskip
Now let $\s$ be a section of $\E$.  Throughout the paper it is convenient
to abbreviate: 
$$
F\,=\,\tfrac12\'|\'\s\'|^2
\tag2-3
$$
If $\{E_i\}$ is a local orthonormal tangent frame in $M$ then by the
defining properties \cite{11} of the connection map, and
holonomy-invariance of $\<\,,\>$ we have:
$$
\align
|\'\dvs\'|^2\,
&=\,\sum_ih\bigl(\dvs(E_i),\dvs(E_i)\bigr) \\
&=\,\sum_iw^p(\s)\left(\bigl\<K\circ d\s(E_i),K\circ d\s(E_i)\bigr\>\,
+\vphantom{\nab{E_i}Y}\,q\'\bigl\<K\circ d\s(E_i),\s\bigr\>^2\right) \\
&=\,w^p(\s)\sum_i\left(\bigl\<\nab{E_i}\s,\nab{E_i}\s\bigr\>\,
+\,q\'\bigl\<\nab{E_i}\s,\s\bigr\>^2\right)  \\
&=\,w^p(\s)\left(|\'\nabla\s\'|^2\,
+\,q\,|\'\nabla F|^2\right),
\tag2-4
\endalign
$$
where $\nabla F$ is the gradient vector.  In the Sasaki case (2-4) reduces
to:
$$
|\'\dvs\'|^2\,=\,|\'\nabla\s\'|^2,
\tag2-5
$$
whereas in the Cheeger-Gromoll case:
$$
|\'\dvs\'|^2\,
=\,\frac{1}{1+|\'\s\'|^2}
\left(|\'\nabla\s\'|^2\,
+\,|\'\nabla F|^2\right)
$$
If $q<0$ and $|\'\s(x)|^2<-1/q$ for all $x\in M$ then the fact that
$h_{p,q}$ is a Riemannian metric in $B\E(1/\sqrt{-q}\,)$ (so that
$|\'\dvs\'|^2\geqs0$), and $d^v\s=0$ if and only if $\nabla\s=0$, allows us
to immediately  deduce from (2-4) {\sl Kato's inequality:}
$$
|\'\nabla\s\'|^2\,+\,q\,|\'\nabla F|^2\,\geqs\,0,
\quad\text{with equality if and only if $\nabla\s=0$.}
\tag2-6
$$
It is not hard to see that (2-6) remains true if $|\'\s\'|^2\leqs-1/q$.
\subhead
Definitions
\endsubhead
For any $q\in\Bbb R$, a smooth section $\s$ satisfying
$\,q\'|\'\s(x)|^2\geqs-1\,$ for all $x\in M$ will be called 
{\sl $q$-Riemannian,} and the set of all such $\s$ will be denoted
$\C(\E,q)$.   Smooth sections satisfying (2-6) will be called 
{\sl $q$-positive,} and the set of all such  will be denoted by $\C(\E,q+)$.
\subhead
Remark 2.4
\endsubhead
If $q\geqs0$ then $\C(\E,q)=\C(\E)$, the space of all smooth sections of
$\E$, and if $q<0$ then $\s$ is $q$-Riemannian precisely when its image lies
in the closure of  the ``Riemannian tube'' $B\E(1/\sqrt{-q}\,)$.  Although
there are tangent vectors (to $\E$) at points on the boundary
$S\E(1/\sqrt{-q}\,)$ of $B\E(1/\sqrt{-q}\,)$ which are $h_{p,q}$-null,
since a $q$-Riemannian section which encounters $S\E(1/\sqrt{-q}\,)$
does so tangentially it follows from Remark 2.3 that the restriction of
$h_{p,q}$ to $\s(M)$ is indeed a Riemannian metric.  
\subhead
Remark 2.5
\endsubhead
If $q_1<q_2$ then $\C(\E,q_1)\subset\C(\E,q_2)$.  Certainly
$\C(\E,q)\subset\C(\E,q+)$, but $\C(\E,q+)$ also includes (for example)
all sections of constant length. 
\medskip
We denote by $E^v_{p,q}$ the vertical energy functional with respect to
$h_{p,q}$.  By (1-4) and (2-4):
$$
E^v_{p,q}(\s)\,
=\,\frac12\int_Mw^p(\s)\left(|\'\nabla\s\'|^2\,
+\,q\,|\'\nabla F|^2\right)\volume(g),
\tag2-7
$$
for all $\s\in\C(\E)$.  We refer to $E^v_{p,q}(\s)$ as the
{\sl vertical $(p,q)$-energy\/} of $\s$.  When $(p,q)$ are understood we
simply write $E^v(\s)$.
\subhead
Remark
\endsubhead
Since $(\E,h_{p,q})\to(M,g)$ is horizontally isometric for all $(p,q)$,
$E^v(\s)$ differs from the full energy $E(\s)$ \cite{10} by a positive
constant, depending only on the dimension and volume of $M$.
\bigskip
\head 
3. The first variation formula
\endhead 
Let $\Sigma\colon M\times\Bbb R\to\E;\Sigma(x,t)=\s_t(x)$ be a smooth
variation of $\s=\s_0$ through sections of $\E$.  Then it is natural to
identify the variation field, which is necessarily vertical, with a family
of sections $\rho_t$ of $\E$:
$$
\rho_t(x)\,
=\,K\circ\frac{d}{dt}\bigl(\s_t(x)\bigr)\,
\in\,\pi^{-1}(x),
$$
bearing in mind that the restriction of $K$ to the vertical distribution is
canonical.  Furthermore, $\Sigma$ may be viewed as a section of the pullback
vector bundle $\pi_1^{-1}\E\to M\times\Bbb R$, where 
$\pi_1\colon M\times\Bbb R\to M;(x,t)\mapsto x$, and this allows the
variation field to be expressed in terms of the pullback connection:
\proclaim{Lemma 3.1}
If $\partial_t$ is the unit vector field on $M\times\Bbb R$ in the positive 
$\Bbb R$-direction then $\nab{\partial_t}\Sigma=\pi_1^{-1}\rho_t$
\rom(the
$\pi_1$-pullback of $\rho_t$\rom), for all $t$. 
\endproclaim
\demo{Proof}
In general, if $f\colon P\to M$ is any smooth map and
$(p,e)\in f^{-1}\E\subset P\times\E$ (thus $f(p)=\pi(e)$) then the
tangent space of $f^{-1}\E$ is the following subspace of 
$T_{(p,e)}(P\times\E)$:
$$
T_{(p,e)}(f^{-1}\E)\,
=\,\{(Y,A):Y\in T_pP,\;A\in T_e\E,\;
df(Y)=d\pi(A)\},
$$
using the natural identification of $T_{(p,e)}(P\times\E)$ with
$T_pP\oplus T_e\E$, and the connection map for the pullback connection is:
$$
\tilde K(Y,A)\,=\,(p,KA)
$$
Therefore, since as a section of $\pi_1^{-1}\E$ we have
$\Sigma(x,t)=\bigl((x,t),\s_t(x)\bigr)$, it follows that:
$$
\nab{\partial_t}\Sigma(x,t)\,
=\,\tilde K\bigl(d\Sigma(\partial_t)\bigr)\,
=\,\bigl((x,t),\,K\circ\frac{d}{dt}\bigl[\'\s_t(x)\'\bigr]\bigr)\,
=\,\bigl((x,t),\rho_t(x)\bigr)
$$
Thus $\nab{\partial_t}\Sigma$ is the $\pi_1$-pullback of $\rho_t$.
\qed
\enddemo
We note also the following property of the curvature of the pullback
connection (in a slightly more general setting):
\proclaim{Lemma 3.2}
Let $\pi_1\colon M\times N\to M;(x,y)\mapsto x$ be the projection from any
product with $M$, and let $\tilde\E=\pi_1^{-1}\E\to M\times N$.  The
curvature of the pullback connection satisfies: 
$$
\tilde R(X,Y)=0,\quad\text{for all $X\in T_xM$, $Y\in T_yN$,}
$$ 
regarding $\,T_xM,T_yN\subset T_{(x,y)}(M\times N)$ in the natural way.
\endproclaim
\demo{Proof}
Extend $X,Y$ to local vector fields in $M,N$ respectively; then $X,Y$ may
also be regarded as local vector fields in $M\times N$.  Suppose
$\tilde e\in\tilde\E_{(x,y)}$; thus $\tilde e=((x,y),e)$ where $e\in\E_x$. 
Extend $e$ to a local section $\a$ of $\E$, and then extend $\tilde e$ to
the local section $\tilde\a(x,y)=((x,y),\a(x))$ of $\tilde\E$.  Thus
$\tilde\a$ is the $\pi_1$-pullback of $\a$, and so by a fundamental
characterization of pullback connections:
$$
\gather
\nab X\tilde\a\,
=\,\bigl((x,y),\nab{d\pi_1(X)}\a\bigr)\,
=\,\widetilde{\nab X\a}, 
\tag3-1 \\
\vspace{1ex}
\nab Y\tilde\a\,
=\,\bigl((x,y),\nab{d\pi_1(Y)}\a\bigr)\,=\,0
\tag3-2
\endgather
$$
Therefore, using the fact that the local vector field $X$ in $M\times N$ is
$\pi_1$-adapted to its counterpart in $M$, and $Y$ is $\pi_1$-adapted to
$0$:
$$
\align
\tilde R(X,Y)\tilde e\,
&=\,\nab X\nab Y\tilde\a-\nab Y\nab X\tilde\a
-\nab{[\'X,Y\']}\tilde\a \\
\vspace{1ex}
&=\,\nab X\nab Y\tilde\a-\nab Y\nab X\tilde\a,
\quad\text{since $[\'X,Y\']=0$,} \\
\vspace{1ex}
&=\,-\nab Y\widetilde{\nab X\a},
\quad\text{by (3-1) and (3-2),} \\
\vspace{1ex}
&=\,0,
\quad\text{by (3-2).}
\tag"\qed"
\endalign
$$
\enddemo
\bigskip
Now it follows from (2-7) that:
$$
\align
\left.\frac{d}{dt}\right|_{t=0}E^v_{p.q}(\s_t)\,
&=\,\frac12\int_M\left.\frac{d}{dt}\right|_{t=0}w^p(\s_t)
\left(|\'\nabla\s\'|^2\,
+\,q\,|\'\nabla F|^2\right)\volume(g)  \\
\vspace{1ex}
&\qquad+\,
\frac12\int_Mw^p(\s)
\left.\frac{d}{dt}\right|_{t=0}\left(|\'\nabla\s_t\'|^2\,
+\,q\,|\'\nabla F_t\'|^2\right)\volume(g)  \\
\vspace{1.5ex}
&=\,V_1+V_2,\quad\text{say.}
\endalign
$$
We abbreviate the variation field $\rho_0=\rho$.  It is also convenient to
define an $\E$-valued 1-form $\v$ on $M$ as follows: 
$$
\v(Y)\,=\,\<\nabla F,Y\>\s\,
=\,\tfrac12\<\nabla|\'\s\'|^2,Y\>\s
\tag3-3
$$
\proclaim{Lemma 3.3}
\medskip\noindent
{\rm(i)}\quad
$\left.\frac{d}{dt}\right|_{t=0}w^p(\s_t)\,
=\,-2p\'w^{p+1}(\s)\'\<\s,\rho\>$
\medskip\noindent
{\rm(ii)}\quad
$\left.\frac{d}{dt}\right|_{t=0}|\'\nabla\s_t\'|^2\,
=\,2\<\nabla\rho,\nabla\s\>$
\medskip\noindent
{\rm(iii)}\quad
$\left.\frac{d}{dt}\right|_{t=0}
|\'\nabla F_t\'|^2\,
=\,2\'\<\v,\nabla\rho\>\,
+\,2\'\bigl\<\nab{\nabla F}\s,\,\rho\bigr\>$
\endproclaim
\demo{Proof}
\flushpar
(i)\quad
We  have:
$$
\align
\left.\frac{d}{dt}\right|_{t=0}w^p(\s_t)\,
&=\,p\'w^{p-1}(\s)\left.\frac{d}{dt}\right|_{t=0}w(\s_t)
=\,p\'w^{p-1}(\s)\,\frac{-1}{(1+|\'\s\'|^2)^2}
\left.\frac{d}{dt}\right|_{t=0}|\'\s_t\'|^2 \\
\vspace{1.5ex}
&=\,-p\'w^{p+1}(\s)\'\left.\frac{d}{dt}\right|_{t=0}|\'\s_t\'|^2
\endalign
$$
By Lemma 3.1:
$$
\frac{d}{dt}\,|\'\s_t\'|^2\,
=\,\frac{d}{dt}\,|\'\Sigma\'|^2\,
=\,2\'\bigl\<\nab{\partial_t}\Sigma,\Sigma\bigr\>\,
=\,2\'\<\'\rho_t,\s_t\>
\tag3-4
$$
\bigskip\noindent
(ii)\quad
Summing over $i$:
$$
\allowdisplaybreaks
\align
\frac12\frac{d}{dt}\,|\'\nabla\s_t\'|^2\,
&=\,\frac12\,\frac{d}{dt}
\bigl\<\nab{E_i}\s_t,\nab{E_i}\s_t\bigr\>\,
=\,\frac12\,\frac{d}{dt}
\bigl\<\nab{E_i}\Sigma,\nab{E_i}\Sigma\bigr\>\,
=\,\bigl\<\nab{\partial_t}\nab{E_i}\Sigma,\nab{E_i}\Sigma\bigr\> \\
\vspace{1ex}
&=\,\bigl\<\nab{E_i}\nab{\partial_t}\Sigma\,
+\,\nab{[\'\partial_t,E_i\']}\Sigma,\,\nab{E_i}\Sigma\bigr\>,
\quad\text{by Lemma 3.2,} \\
\vspace{1ex}
&=\,\bigl\<\nab{E_i}\nab{\partial_t}\Sigma,\nab{E_i}\Sigma\bigr\>,
\quad\text{since $[\'\partial_t,E_i\']=0$,} \\
\vspace{1ex}
&=\,\bigl\<\nab{E_i}\rho_t,\nab{E_i}\s_t\bigr\>,
\quad\text{by Lemma 3.1,} \\
\vspace{1ex}
&=\,\<\nabla\rho_t,\nabla\s_t\>
\endalign
$$
\medskip\noindent
(iii)\quad
Summing over $i$:
$$
\allowdisplaybreaks
\align
\frac12\left.\frac{d}{dt}\right|_{t=0}
|\'\nabla F_t\'|^2\,
&=\,\frac12\left.\frac{d}{dt}\right|_{t=0}
\bigl\<\nab{E_i}\s_t,\,\s_t\bigr\>^2\,
=\,\bigl\<\nab{E_i}\s,\s\bigr\>\left.\frac{d}{dt}\right|_{t=0}
\bigl\<\nab{E_i}\s_t,\s_t\bigr\> \\
\vspace{1ex}
&=\,\bigl\<\nab{E_i}\s,\s\bigr\>\left(
\bigl\<\nab{\partial_t}\nab{E_i}\Sigma,\s\bigr\>\,
+\,\bigl\<\nab{E_i}\s,\nab{\partial_t}\Sigma\bigr\>\right)
\vspace{1ex}
&=\,\frac12(E_i.|\'\s\'|^2)
\left(\bigl\<\nab{E_i}\nab{\partial_t}\Sigma,\s\bigr\>\,
+\,\bigl\<\nab{E_i}\s,\nab{\partial_t}\Sigma\bigr\>\right),
\quad\text{by Lemma 3.2,} \\
\vspace{1ex}
&=\,\<\nabla F,E_i\>
\left(\bigl\<\nab{E_i}\rho,\s\bigr\>\,
+\,\bigl\<\nab{E_i}\s,\rho\bigr\>\right),
\quad\text{by Lemma 3.1,} \\
\vspace{1ex}
&=\,\bigl\<\nab{E_i}\rho,\,\v(E_i)\bigr\>\,
+\,\bigl\<\nab{\nabla F}\s,\,\rho\bigr\>,
\quad\text{by (3-3)} \\
\vspace{1ex}
&=\,\<\v,\nabla\rho\>\,
+\,\bigl\<\nab{\nabla F}\s,\,\rho\bigr\>.
\tag"\qed"
\endalign
$$
\enddemo
\bigskip
It follows from Lemma 3.3 that:
$$
\gather
V_1\,=\,-p\int_Mw^{p+1}(\s)\,
\bigl\<\bigl(|\'\nabla\s\'|^2\,
+\,q\,|\'\nabla F|^2\bigr)\s,
\,\rho\bigr\>\volume(g)
\vphantom{\left(\frac12\right)}
\tag3-5 \\
\vspace{1.5ex}
V_2\,=\,\int_Mw^p(\s)\left(\vphantom{\nab{E_i}}\right.
\bigl\<\nabla\s+q\'\v,\nabla\rho\bigr\>\,
+\,\bigl\<q\,\nab{\nabla F}\s,\,\rho\bigr\>
\left.\vphantom{\nab{E_i}}\right)\volume(g)\vphantom{\left(\frac12\right)}
\tag3-6
\endgather
$$
\medskip\noindent
The expression (3-6) is only partially in divergence form, a situation which
we rectify with the following sequence of calculations.  For any $\E$-valued
1-form $\b$ on $M$ we have (summing over $i$):
$$
\align
\nabla^*(f\b)\,
&=\,-\nab{E_i}(f\b)(E_i)\,
=\,-df(E_i)\,\b(E_i)\,
-\,f\,\nab{E_i}\b(E_i) \\
\vspace{1ex}
&=\,f\,\nabla^*\b\,
-\,\b(\nabla f)
\endalign
$$
In particular, if $f=w^p(\s)$ then: 
$$
\align
df\,&=\,p\'w^{p-1}(\s)\,d\bigl(w(\s)\bigr)\,
=\,p\'w^{p-1}(\s)\,\frac{-2}{(1+|\'\s\'|^2)^2}\,\<\nabla\s,\s\> \\
&=\,-2p\'w^{p+1}(\s)\,\<\nabla\s,\s\>,
\endalign
$$ 
and so, summing over $i$:
$$
\align
\nabla f\,&=\,df(E_i)E_i\,
=\,-2p\'w^{p+1}(\s)\'\bigl\<\nab{E_i}\s,\s\bigr\>E_i\,
=\,-p\'w^{p+1}(\s)\'(E_i.|\'\s\'|^2)\'E_i \\
\vspace{1ex}
&=\,-p\'w^{p+1}(\s)\'\nabla|\'\s\'|^2\,
=\,-2p\'w^{p+1}(\s)\'\nabla F
\endalign
$$
Taking $\b=\nabla\s+q\'\v$ yields:
$$
\nabla^*(f\b)\,
=\,w^p(\s)\bigl(\nabla^*\nabla\s\,
+\,q\,\nabla^*\v\bigr)\,
+\,2p\'w^{p+1}(\s)\left(\nab{\nabla F}\s\,
+\,q\'\v\bigl(\nabla F\bigr)\right)
\tag3-7
$$
\proclaim{Lemma 3.4}
\quad
$\<\nabla^*\nabla\s,\s\>\,
=\,|\'\nabla\s\'|^2\,+\,\Delta F$
\endproclaim
\demo{Proof}
Abbreviating $|\'\s\'|^2=2F$:
$$
\align
\<\nabla^*\nabla\s,\s\>\,
&=\,-\sum_i\bigl\<\nab{E_i}\nab{E_i}\s\,
-\,\nab{\nab{E_i}E_i}\s,\,\s\bigr\> \\
&=\,-\sum_i\left(E_i.\bigl\<\nab{E_i}\s,\s\bigr\>\,
-\,\bigl\<\nab{E_i}\s,\nab{E_i}\s\bigr\>\,
-\,\bigl\<\nab{\nab{E_i}E_i}\s,\,\s\bigr\>\right) \\
&=\,|\'\nabla\s\'|^2\,
-\,\sum_i\left(E_i.E_i.F\,
-\,(\nab{E_i}E_i).F\right) \\
&=\,|\'\nabla\s\'|^2\,
-\,\Tr\nabla dF\,
=\,|\'\nabla\s\'|^2\,
+\,\Delta F.
\tag"\qed"
\endalign
$$
\enddemo
\subhead
Note
\endsubhead
Our sign convention for the Laplace-Beltrami operator is:
$$
\Delta F\,=\,-\Tr\nabla dF
$$
\proclaim{Lemma 3.5}
\quad $\nabla^*\v\,=\,\bigl(\Delta F\bigr)\s\,-\,\nab{\nabla F}\s$,\quad
where $2F=|\'\s\'|^2$.
\endproclaim
\demo{Proof}
Summing over $i$:
$$
\allowdisplaybreaks
\align
\nabla^*\v\,
&=\,-\nab{E_i}\v(E_i)\,
=\,-\nab{E_i}(\v(E_i)),
\quad\text{if $\{E_i\}$ is at the centre of a normal chart,} \\
\vspace{1ex}
&=\,-\nab{E_i}\bigl(\<\nabla F,E_i\>\s\bigr)\,
=\,-\nab{E_i}\bigl(\<\nabla F,E_i\>\s\bigr)\,
=\,-\nab{E_i}\bigl((E_i.F)\s\bigr) \\
\vspace{1ex}
&=\,-\nab{E_i}\left(\bigl\<\nab{E_i}\s,\s\bigr\>\s\right)\,
=\,\bigl\<-\nab{E_i}\nab{E_i}\s,\s\bigr\>\'\s\,
-\,|\'\nabla\s\'|^2\'\s\,
-\,\bigl\<\nab{E_i}\s,\s\bigr\>\'\nab{E_i}\s \\
\vspace{1ex}
&=\,\bigl\<\nabla^*\nabla\s,\s\bigr\>\'\s\,
-\,|\'\nabla\s\'|^2\'\s\,
-\,\nab{\nabla F}\s \\
\vspace{1ex}
&=\,\bigl(\Delta F\bigr)\s\,
-\,\nab{\nabla F}\s,
\quad\text{by (2-3) and Lemma 3.4.}
\tag"\qed"
\endalign
$$
\enddemo
\proclaim{Theorem 3.6}
For any $1$-parameter smooth variation $\s_t$ of $\s$ through sections we
have:
$$
\aligned
\left.\frac{d}{dt}\right|_{t=0}E^v_{p,q}(\s_t)\,
&=\,\int_M\left(w^p(\s)\bigl\<\nabla^*\nabla\s\,
+\,q\,\Delta F\,\s,\,\rho\bigr\>
\right. \\
&\qquad\quad
\left.+\,\,p\'w^{p+1}(\s)\,\bigl\<2\,\nab{\nabla F}\s\,
+\,q\,|\'\nabla F|^2\'\s\,
-\,|\'\nabla\s\'|^2\'\s,\,\rho\bigr\>
\right)\volume(g),
\endaligned
$$
where $2F=|\'\s\'|^2$.  For all $(p,q)$, $\s$ is a $(p,q)$-harmonic section
of $\E$ if and only if:
$$
T_p(\s)\,=\,\phi_{p,q}(\s)\'\s,
$$
where:
$$
\align
T_p(\s)\,&=\,(1+2F)\'\nabla^*\nabla\s\,
+\,2p\,\nab{\nabla F}\s, \\
\vspace{1ex}
\phi_{p,q}(\s)\,
&=\,p\,|\'\nabla\s\'|^2\,-\,pq\,|\'\nabla F|^2\,
-\,q(1+ 2F)\Delta F.
\endalign
$$
\endproclaim
\demo{Proof}
\noindent
Applying Lemma 3.5 to equation (3-7), and noting the cancellation of terms
involving $\;\pm q\,\nab{\nabla F}\s$, yields:
$$
\aligned
V_2\,&=\,\int_M\left(\vphantom{\nab{\nabla F}\s}
w^p(\s)\,\bigl\<\nabla^*\nabla\s\,
+\,q\'(\Delta F)\'\s,\,\rho\bigr\>\right. \\
&\qquad\qquad\qquad
\left.+\,2p\'w^{p+1}(\s)\,\bigl\<\nab{\nabla F}\s\,
+\,q\,|\'\nabla F|^2\'\s,\,\rho\bigr\>
\right)\volume(g)
\endaligned
$$
The first variation formula is then obtained by combining this
with (3-5).
\qed
\enddemo
\subhead
Remark 3.7
\endsubhead
By (2-6) and (2-7), $E^v(\s)\geqs0$ for all $q$-positive sections $\s$, and
the zeroes of $E^v$ in $\C(\E,q+)$ are precisely the parallel sections,
which are therefore the absolute minima of the restriction of $E^v$ to
$\C(\E,q+)$.  In particular, any parallel $\s$ is a $(p,q)$-harmonic
section of $\E$ for all $(p,q)$ with $q\'|\'\s\'|^2>-1$.  However, it
follows from Theorem 3.6 that parallel sections are in fact
$(p,q)$-harmonic for all $(p,q)$.
\subhead
Remark
\endsubhead
In the Sasaki case, when $(p,q)=(0,0)$, we recover the Euler-Lagrange
equations (1-2).  In fact we get:
$$
\left.\frac{d}{dt}\right|_{t=0}E^v(\s_t)\,
=\,\int_M\<\nabla^*\nabla\s,\rho\>\volume(g)
\tag3-8
$$ 
\subhead
Remark 3.8
\endsubhead
If $|\'\s\'|=k$ (constant) then $\nabla F=0$ and it follows from (2-7) that
$E^v(\s)$ is a constant multiple (depending only on $k$ and $p$) of the
Sasaki vertical energy.  Therefore if $k>0$ then $\s$ is a
$(p,q)$-harmonic section of $S\E(k)$ (ie\. a critical point of $E^v$ with
respect to $h_{p,q}$ and smooth variations through sections of $S\E(k)$) if
and only if $\s$ is a harmonic section of $S\E(k)$ with respect to the
Sasaki metric, so we simply say that $\s$ is a {\sl harmonic section of
$S\E(k)$.}  Differentiating the constraint equation $|\'\s_t\'|^2=k^2$ with
respect to $t$ yields $\<\s,\rho\>=0$, by (3-4), and it therefore follows
from (3-8) that $\s$ is a harmonic section of $S\E(k)$ if and only if
$\,\nabla^*\nabla\s=\lambda\,\s$ for some smooth 
$\lambda\colon M\to\Bbb R$.  It then follows from Lemma 3.4 that:
$$
k^2\lambda\,
=\,\<\nabla^*\nabla\s,\s\>\,
=\,|\'\nabla\s\'|^2
$$  
Thus $\s$ is a harmonic section of $S\E(k)$ if and only if $\s$ satisfies
equation (1-3).
\bigskip
\head 
4. Main Theorems
\endhead 
Our first result is Theorem A of the Introduction.
\proclaim{Theorem 4.1}
Suppose $\s$ has constant length $k>0$.  Then $\s$ is a $(p,q)$-harmonic
section of $\E$ if and only if $\s$ is parallel, except when
$\,p=1+1/k^2\,$ in which case $\s$ is a $(p,q)$-harmonic section of
$\E$ if and only if $\s$ is a harmonic section of $S\E(k)$.
\endproclaim
\demo{Proof}
If $|\'\s\'|=k$ then:
$$
T_p(\s)\,=\,(1+k^2)\'\nabla^*\nabla\s
\quad\text{and}\quad
\phi_{p,q}(\s)\,=\,p\'|\'\nabla\s\'|^2
$$
Hence the $(p,q)$-harmonic section equations reduce to:
$$
\nabla^*\nabla\s\,
=\,\frac{p}{1+k^2}\,|\'\nabla\s\'|^2\'\s,
$$
which by Lemma 3.4 implies:
$$
\frac{(p-1)k^2-1}{k^2+1}\,|\'\nabla\s\'|^2\,=\,0
$$
Therefore $\,\nabla\s=0$, except when $\,p=1+1/k^2$, in which case the
$(p,q)$-harmonic section equations become:
$$
\nabla^*\nabla\s\,=\,(p-1)|\'\nabla\s\'|^2\s
\tag4-1
$$
But (4-1) is precisely the equation (1-3) for $\s$ to be a harmonic section
of $S\E(k)$ when $k=1/\sqrt{p-1}\'$.
\qed
\enddemo
\subhead
Example 4.2
\endsubhead
If $\xi$ is a harmonic section of the unit sphere bundle $S\E(1)$ then it
follows from equation (1-3) that $\s=k\'\xi$ is a harmonic section of
$S\E(k)$ for all $k>0$.  It follows from Theorem 4.1 that $\s$ is a
$(p,q)$-harmonic section of $\E$ for $p=1+1/k^2$ and all $q$.  Thus given
a harmonic section of $S\E(1)$ and any $(p,q)$ with $p>1$ it is possible to
construct a $(p,q)$-harmonic section of $\E$.  For example, let
$M=S^{2m+1}$ and let $\xi$ be the standard Hopf vector field: $\xi(x)=ix$
(where $i=\sqrt{-1}$, thinking of  
$x\in S^{2m+1}\subset\Bbb R^{2m+2}\cong\Bbb C^{m+1}$).  Then $\xi$ is a
harmonic section of $S\E(1)$ where $\E=TM$ \cite{22}, and so
$\s=\xi/\sqrt{p-1}\,$ is a $(p,q)$-harmonic section of $\E$.  In
particular, this shows that for all $(p,q)$ with $p>1$ there exist examples
of $(p,q)$-harmonic sections of constant length which are not parallel. 
\medskip
For non-compact $M$, sections of constant length are generalized by those
for which $|\'\s\'|^2\colon M\to\Bbb R$ is a harmonic function.  Lemma 3.4
shows that if such a section is $(0,0)$-harmonic then it is parallel.  We
now investigate this ``Bernstein phenomenon'' for other values of $(p,q)$. 
For this it is convenient to identify the following regions of the
$(p,q)$-plane, for any $\mu>0$:
$$
\gather
\F_-(\mu)\,=\,\{(p,q):p<0,\,\mu q\leqs2p\},\qquad
\F_0\,=\,\{(p,q):0\leqs p\leqs1\}, \\
\vspace{1ex}
\F_1(\mu)\,=\,\{(p,q):p>1,\,\mu q<1-p\},
\endgather
$$
and then define:
$$
\F(\mu)\,=\,\F_-(\mu)\,\cup\,\F_0\,\cup\,\F_1(\mu)
$$
\proclaim{Theorem 4.3}
Suppose that $\mu\geqs1/2$ and $(p,q)\in\F(\mu)$, $\s$ is a $\mu
q$-Riemannian section of $\E$, and $|\'\s\'|^2$ is a harmonic function. 
Then $\s$ is a $(p,q)$-harmonic section if and only if $\s$ is parallel.
\endproclaim
\demo{Proof}
From Theorem 3.6 and Lemma 3.4:
$$
\align
\<T_p(\s),\s\>\,
&=\,(1+2F)\left(|\'\nabla\s\'|^2\,
+\,\Delta F\right)\,
+\,2p\,|\'\nabla F|^2 \\
\<\phi_{p,q}(\s)\'\s,\s\>\,
&=\,2pF\,|\'\nabla\s\'|^2\,
-\,2pqF\,|\'\nabla F|^2\,
-\,2qF(1+2F)\Delta F
\endalign
$$
Therefore:
$$
\align
\<T_p(\s)-\phi_{p,q}(\s)\'\s,\s\>\,
&=\,\bigl(1+2(1-p)F\bigr)|\'\nabla\s\'|^2\,
+\,(1+2qF)(1+2F)\Delta F \\
&\qquad
+\,2p(1+qF)\'|\'\nabla F|^2 \\
\vspace{1ex}
&=\,C_1\,|\'\nabla\s\'|^2\,
+\,C_2\,\Delta F\,
+\,C_3\,|\'\nabla F|^2,
\quad\text{say.}
\tag4-2
\endalign
$$
Thus if $\s$ is $(p,q)$-harmonic and $F$ is a harmonic function then by
Theorem 3.6:
$$
0\,=\,C_1\,|\'\nabla\s\'|^2\,
+\,C_3\,|\'\nabla F|^2
\tag4-3
$$
We now consider each subregion of $\F(\mu)$ separately.
\flushpar
(i)\quad
$(p,q)\in\F_{-}(\mu)$
\par
Since $C_1>0$ and $\s$ is $\mu q$-Riemannian, applying Kato's inequality
(2-6) to (4-3) yields:
$$
\align
0\,&\geqs\,(C_3-\mu q\'C_1)\,|\'\nabla F|^2 \\
&=\,\bigl(2p-\mu q\,
+\,2((\mu+1)p-\mu)qF\bigr)\'|\'\nabla F|^2 \\
&=\,(A\,+\,B\,qF)\'|\'\nabla F|^2,
\quad\text{say.}
\tag4-4
\endalign
$$
We have $A\geqs0$ and:
$$
B\,qF\,=\,2\bigl(\mu(p-1)+p\bigr)qF\,>\,0,
$$
since $\mu>0$ and $p,q<0$.
Hence
$\nabla F=0$ identically, and so $F$ is constant.  But then $\s$ is parallel,
by Theorem 4.1.
\smallskip\noindent
(ii)\quad
$(p,q)\in\F_0$
\par
We have $C_1>0$ in (4-3).
If $q\geqs0$ then $C_3\geqs0$, and therefore $\nabla\s=0$.  If $q<0$ then
since $C_1>0$ the inequality (4-4) still holds.  For $p\leqs\mu/(\mu+1)$ we
have $B\,qF\geqs0$, and since $A>0$ it follows as in (i) that $\s$
is parallel.  If $p>\mu/(\mu+1)$ then $B>0$ and since $\s$ is 
$\mu q$-Riemannian we have:
$$
B\,qF\,\geqs\, B\left(-\frac{1}{2\mu}\right)\,
=\,1\,-\,\left(\frac{\mu+1}{\mu}\right)p
$$
Therefore (4-4) may be strengthened:
$$
0\,\geqs
\left(2p-\mu q+1-\left(\frac{\mu+1}{\mu}\right)p\right)|\'\nabla F|^2\,
=\,\left(\frac{\mu(p+1)-p}{\mu}\,
-\,\mu q\right)|\'\nabla F|^2
$$
Since $\mu(p+1)-p\geqs0$ for all $p\in[\'0,1\']$ if and only if
$\mu\geqs1/2$, the coefficient of $|\'\nabla F|^2$ is strictly positive.  It
therefore follows as in (i) that $\s$ is parallel.
\smallskip\noindent
(iii)\quad
$(p,q)\in\F_1(\mu)$
\par
Since $\mu q<1-p$ and $\s$ is $\mu q$-Riemannian, it follows that $C_1>0$
in (4-3).  Furthermore since $\mu\geqs1/2$ we have $\mu q\leqs q/2$,
hence $\s$ is $(q/2)$-Riemannian.  Therefore $C_3\geqs0$ in (4-3). Hence
$\nabla\s=0$.
\qed
\enddemo
\subhead
Remark
\endsubhead
Of the $(p,q)$-harmonic sections addressed by Theorem 4.3, only those with
$\mu\geqs1$ are necessarily $q$-Riemannian.  Thus the ``Bernstein
phenomenon'' described by Theorem 4.3 is not necessarily confined to
$(p,q)$-harmonic sections whose image lies in the region of $\E$ where
$h_{p,q}$ is a Riemannian metric.
\subhead
Remark
\endsubhead
If $|\'\s\'|=k>0$ and $\s$ is $\mu q$-Riemannian then $1/k^2\geqs-\mu q$. 
Hence if in addition $p=1+1/k^2$ then $1-p\leqs\mu q$.  It therefore follows
from Theorem 4.1 that $\F_1(\mu)$ is the best possible ``Bernstein region''
when $p>1$ in Theorem 4.3.
\subhead
Remark 4.4
\endsubhead
It follows from Theorem 4.3 that if  $\s$ is a harmonic section with respect
to the Riemannian metric $h_{p,p}$ for any $0\leqs p\leqs1$ and $|\'\s\'|^2$
is a harmonic function, then $\s$ is parallel.  In particular, this is the
case for the Cheeger-Gromoll metric ($p=1$).
\medskip
Remark 4.4 illustrates a surprising feature of Theorem 4.3, which is that
for all $(p,q)$ in the following vertical strip:
$$
\F_0^{+}\,=\,\{(p,q):0\leqs p\leqs1,\,q\geqs0\},
$$
any $(p,q)$-harmonic section $\s$ with $|\'\s\'|^2$ harmonic is necessarily
parallel, without any {\it a priori\/} bound on its length.  By placing an
alternative bound on $|\'\s\'|^2$ when $q\geqs0$ and $p>1$ it is possible to
extend $\F_0^{+}$ rightwards into the following adjacent region:
$$
\G_1(\nu)\,=\,\{(p,q):p>1,\,q\geqs2\nu(1-p)\},
$$
for any $\nu>0$.
\subhead
Definition
\endsubhead
A $q$-Riemannian section $\s$ is {\sl strictly $q$-Riemannian\/} if
$\,q\'|\'\s(x)|^2>-1$ for some $x\in M$.
\subhead
Remark
\endsubhead
For this condition to have any force it is clearly necessary for $M$ to be
connected.  In our next two results (Theorems 4.5\,(b) and 4.6\,(b))
connectedness of $M$ is therefore essential.
\proclaim{Theorem 4.5}
Suppose $(p,q)\in\G_1(\nu)$ and either:
\flushpar
\rom{(a)}\quad
$\nu>1$ and $\s$ is a $\nu(1-p)$-Riemannian section of $\E$\rom;
\flushpar
\rom{(b)}\quad
$\nu=1$ and $\s$ is a strictly $(1-p)$-Riemannian section of $\E$\rom.
\flushpar
Suppose also that $|\'\s\'|^2$ is a harmonic function on $M$.  Then $\s$ is
a $(p,q)$-harmonic section of $\E$ if and only if $\s$ is parallel.
\endproclaim
\demo{Proof}
Since  $2\nu(1-p)F\geqs-1$ and $q\geqs2\nu(1-p)$ we have $qF\geqs-1$; hence
$C_3$ in (4-3) is non-negative.  Furthermore since $\nu\geqs1$ and $\s$ is
$\nu(1-p)$-Riemannian, $\s$ is $(1-p)$-Riemannian (Remark 2.5); hence $C_1$
in (4-3) is also non-negative.  Now consider $\,U=\{x\in M:C_1(x)>0\}$.  It
follows from (4-3) that $\nabla\s$ vanishes on $U$. If $\nu>1$ we have
$U=M$ and hence $\s$ is parallel.  If $\nu=1$ we have in particular that
$|\'\s\'|$ is constant on $U$, hence $U$ is closed.   But $U$ is open, and
non-empty by hypothesis.  Therefore $U=M$, since $M$ is connected.  Thus
$\s$ is parallel.
\qed
\enddemo
\subhead
Remark
\endsubhead
The $(p,q)$-harmonic sections addressed by Theorem 4.5 are $q$-Riemannian if
$q\geqs\nu(1-p)$, but not necessarily so if $2\nu(1-p)\leqs q<\nu(1-p)$.
\subhead
Remark
\endsubhead
For all $\mu\geqs1/2$ and $\nu\geqs1$ we have:
$$
\F_1(\mu)\,\cup\,\G_1(\nu)\,=\,\{(p,q):p>1\},
$$
and:
$$
\F_1(\mu)\,\cap\,\G_1(\nu)\,
=\,\{(p,q):p>1,\,2\mu\nu(1-p)\leqs\mu q<1-p\}\,
=\,\W(\mu,\nu),\quad\text{say.}
$$
In particular, the union is disjoint precisely when $\mu=1/2$ and $\nu=1$.
Furthermore, for all $(p,q)\in\W(\mu,\nu)$ every $\mu q$-Riemannian section
is strictly $(1-p)$-Riemannian, so that Theorem 4.5 is consistent with
Theorem 4.3 in the region where both theorems apply, namely $\W(\mu,1)$. 
In addition, it follows from Theorem 4.5 that Theorem 4.3 extends to the
closure of $\F(\mu)$ for all sections which are strictly
$\mu q$-Riemannian.
\medskip
We now consider the compact case.  To identify a ``Bernstein region'' of the
$(p,q)$-plane it is convenient to introduce the following monotone
decreasing piecewise-linear cut-off function:
\subhead
Definition
\endsubhead
For any $\nu>0$ define $\varrho_\nu\colon[-4,\infty)\to\Bbb R$ by:
$$
\varrho_\nu(p)\,=\,
\cases
-1-p,&\text{if $\,-4\leqs p\leqs-1$,} \\
0,&\text{if $\,-1\leqs p\leqs2$,} \\
\nu(2-p)/2,&\text{if $\,p\geqs2$.}
\endcases
$$
\proclaim{Theorem 4.6}
Suppose $M$ is compact, $p\geqs-4$, $q\geqs\varrho_\nu(p)$ and either:
\flushpar
\rom{(a)}\quad
$\nu>1$ and $\s$ is a $\nu(1-p)$-Riemannian section of $\E$\rom;
\flushpar
\rom{(b)}\quad
$\nu=1$ and $\s$ is a strictly $(1-p)$-Riemannian section of $\E$\rom.
\flushpar
Then $\s$ is $(p,q)$-harmonic if and only if $\s$ is parallel.
\endproclaim
\remark{Remark}
If $p\leqs1$ then all sections are strictly $(1-p)$-Riemannian.
\endremark
\demo{Proof}
Referring to identity (4-2), we note first that repeated use of the
Divergence Theorem yields the following integral formula for the term
involving the Laplacian (where all integrals are taken over $M$, with
respect to the Riemannian volume element):
$$
\allowdisplaybreaks
\align
\int C_2\,\Delta F\,
&=\,\int(1+2F)\Delta F\,
+\,2q\int F(1+2F)\Delta F \\
\vspace{1ex}
&=\,2\int|\'\nabla F\'|^2\,
+\,2q\int\bigl\{|\'\nabla F\'|^2\,
+\,2g(\nabla(F^2),\nabla F)\bigr\} \\
\vspace{1ex}
&=\,2(1+q)\int|\'\nabla F\'|^2\,
+\,8q\int F\,|\'\nabla F\'|^2 
\endalign
$$
Therefore if $\s$ is $(p,q)$-harmonic then integration of (4-2) yields:
$$
0\,=\,\int C_1\,|\'\nabla\s\'|^2\,
+\,2(p+q+1)\int|\'\nabla F|^2\,
+\,2(p+4)q\int F\,|\'\nabla F|^2
\tag4-5
$$
Since $\s$ is $\nu(1-p)$-Riemannian and $\nu\geqs1$, $\s$ is
$(1-p)$-Riemannian (Remark 2.5); hence $C_1\geqs0$.  Define $\,U=\{x\in
M:C_1(x)>0\}$.  If $p\in[-4,2\']$ then $q\geqs\varrho_\nu(p)$ if and only if
$p+q+1\geqs0$ and $(p+4)q\geqs0$.  It therefore follows from (4-5) that
$C_1\'|\'\nabla\s\'|^2$ vanishes identically, and so $\nabla\s$ vanishes on
$U$.  If $p\geqs2$ then since $\nu\geqs1$:
$$
\nu(1-p)\,<\,\nu\left(1-\frac p2\right)\,
=\,\varrho_\nu(p)\,\leqs\,q
\tag4-6
$$ 
Therefore, since $\s$ is $\nu(1-p)$-Riemannian, $\s$ is also
$q$-Riemannian.  It then follows from (4-5) that:
$$
0\,\geqs\,\int C_1\,|\'\nabla\s\'|^2\,
+\,(p+2q-2)\int|\'\nabla F|^2\,
\geqs\,\int C_1\,|\'\nabla\s\'|^2,
$$
since by (4-6):
$$
2q\,\geqs\,\nu(2-p)\,\geqs\,2-p,
\quad\text{because $\nu\geqs1$.}
$$
We therefore deduce again that $\nabla\s$ vanishes on $U$.  If $\nu>1$
then $U=M$.  If $\nu=1$ then the connectedness argument of Theorem 4.5
may be used to conclude the proof. 
\qed
\enddemo
\subhead
Remark
\endsubhead
All the $(p,q)$-harmonic sections addressed by Theorem 4.6 are in fact
$q$-Riemannian (as shown during the proof).
\subhead
Remark
\endsubhead
The re-scaled Hopf vector fields described in Example 4.2
show that Theorem 4.6 is false if $\s$ is merely $(1-p)$-Riemannian.
\subhead
Remark 4.7
\endsubhead
It follows from Theorem 4.6 that if $0\leqs p\leqs1$ then any
$(p,p)$-harmonic section of a compact vector bundle is parallel.  In
particular, this is true of both the Sasaki ($p=0$) and Cheeger-Gromoll
($p=1$) metrics \cite{16} (cf\. Remark 4.4).
\medskip
It follows from Theorem 4.1 that if a section $\s$ of constant length is
$(p,q)$-harmonic then $\s$ is $(p,r)$-harmonic for all $r\in\Bbb R$. 
However the following result shows that this is exceptional.
\proclaim{Theorem 4.8}
Suppose $M$ is compact, and $\s$ is a section of $\E$ whose length is not
constant.  Then for each $p\in\Bbb R$ there exists at most one $q\in\Bbb R$
such that $\s$ is $(p,q)$-harmonic.
\endproclaim
\demo{Proof}
If $\s$ is both $(p,q)$-harmonic and $(p,r)$-harmonic then it follows from
Theorem 3.6 that:
$$
\align
0\,&=\,\phi_{p,r}(\s)\,
-\,\phi_{p,q}(\s)\,
=\,p(q-r)\'|\'\nabla F|^2\,
+\,(q-r)(1+2F)\Delta F \\
\vspace{1ex}
&=\,(q-r)\left(p\,|\'\nabla F\'|^2\,
+\,(1+2F)\Delta F\right)
\tag4-7
\endalign
$$
Then by the Divergence Theorem:
$$
0\,=\,(q-r)\int\left(p\,|\'\nabla F\'|^2\,
+\,(1+2F)\Delta F\right)\,
=\,(q-r)(p+2)\int|\'\nabla F\'|^2
$$
If $q\neq r$ then since $F$ is not constant it follows that $p=-2$, in
which case (4-7) reduces to:
$$
0\,=\,\Delta F\,+\,2F\Delta F\,-\,2\,|\'\nabla F\'|^2\,
=\,\Delta(F+F^2),
$$
where we have used the formula for the Laplacian of a product:
$$
\Delta(f_1\'f_2)\,
=\,f_1\,\Delta f_2\,+\,f_2\,\Delta f_1\,
-\,2g(\nabla f_1,\nabla f_2)
\tag4-8
$$
But this implies that $F^2+F$, and hence $F$, is constant.
\qed
\enddemo
\subhead
Remark
\endsubhead
It follows directly from (4-7) that Theorem 4.8 also holds if $M$ is
non-compact provided $|\'\s\'|^2$ is a non-constant harmonic function.
\medskip
Combining Theorems 4.6 and 4.8 yields the following result, which may be
regarded as a partial complement to Theorem 4.1:
\proclaim{Corollary 4.9}
Suppose $M$ is compact, and $\s$ is a section of $\E$ whose length is not
constant.  Then for each $p\in\Bbb R$ there exists at most one $q\in\Bbb R$
such that $\s$ is $(p,q)$-harmonic, and if either $-4\leqs p\leqs1$, or
$p>1$ and $\|\'\s\'\|_\infty^2\leqs1/\nu(p-1)$ for some $\nu\geqs1$, then
$q<\varrho_\nu(p)$.
\endproclaim
Theorem B is now a simple consequence of Corollary 4.9 (with $\nu=1$).
Combining Theorem 4.3 (with $\mu=1/2$), Theorem 4.5 (with $\nu=1$) and
Theorem 4.8 yields the following non-compact analogue of Corollary 4.9:
\proclaim{Corollary 4.10}
Suppose $\s$ is a section of $\E$ for which $|\'\s\'|^2$ is a non-constant
harmonic function.  Then for each $p\in\Bbb R$ there exists at most one
$q\in\Bbb R$ such that $\s$ is $(p,q)$-harmonic, and:
\flushpar
\rom{(a)}\quad
if $\,p<0$ and $\,q\leqs4p$ then $|\'\s(x)|^2>-2/q$ for some $x\in M$\rom;
\flushpar
\rom{(b)}\quad
if $\,0\leqs p\leqs1$ then $\,q<0$ and $|\'\s(x)|^2>-2/q$ for some
$x\in M$\rom;
\flushpar
\rom{(c)}\quad
if $\,p>1$ and $\,q<2(1-p)$ then $|\'\s(x)|^2>-2/q$ for some $x\in M$\rom;
\flushpar
\rom{(d)}\quad
if $\,p>1$ and $\,q\geqs2(1-p)$ then $|\'\s(x)|^2>1/(p-1)$ for some
$x\in M$\rom.
\endproclaim
\bigskip
\head 
5. Example: Vector Fields on Spheres
\endhead
Let $M=S^n\subset\Bbb R^{n+1}$ be the unit sphere, with the induced
Riemannian metric $g$, and let $\E=TM$ with $\<\,,\>=g$.  We look for
$(p,q)$-harmonic sections of $\E$ amongst the class of {\sl conformal
gradient fields\/} on $M$.  Thus, let $a\in\Bbb R^{n+1}\smallsetminus\{0\}$,
and let $\lambda\colon S^n\to\Bbb R$ be the restriction of the linear
functional on $\Bbb R^{n+1}$ dual to $a$.  Define $\s=\nabla\lambda$.  
We refer to $a$ as the {\sl axial vector\/} of $\s$, and note that 
$|\'a\'|=\|\'\s\'\|_\infty$.  We say that $\s$ is {\sl standard\/} if
$|\'a\'|=1$.  The following is a collation of data required to compute the
$(p,q)$-harmonic section equations.
\proclaim{Lemma 5.1}
Suppose the axial vector of $\s$ has length $c$.  Then: 
\flushpar
{\rm(i)}\quad
$\nab X\s\,=\,-\lambda\'X$
\flushpar
{\rm(ii)}\quad
$\nabla^*\nabla\s\,=\,\s$
\flushpar
{\rm(iii)}\quad
$2F\,=\,c^2-\lambda^2$
\flushpar
{\rm(iv)}\quad
$\nabla F\,=\,-\lambda\'\s$
\flushpar
{\rm(v)}\quad
$\Delta F\,=\,c^2-(n+1)\lambda^2$
\endproclaim
\demo{Proof}
(i) and (ii) are well-known; see for example \cite{25}.
\flushpar
(iii)\quad
By definition:
$$
\s(x)\,=\,a-\lambda(x)x\,
=\,a-(a\dotprod x)x,
$$ 
where $a\dotprod x$ is the Euclidean dot product, and it therefore follows
that:
$$
2F(x)\,=\,|\'\s(x)|^2\,
=\,|\'a\'|^2-2(a\dotprod x)^2+(a\dotprod x)^2\,
=\,c^2-\lambda(x)^2
$$
(iv)\quad
From (iii):
$$
\nabla F\,=\,\nabla F\,
=\,-\frac12\,\nabla\lambda^2\,
=\,-\lambda\,\nabla\lambda\,=\,-\lambda\'\s
$$
(v)\quad
We have:
$$
\align
\Delta F\,&=\,-\div\nabla F\,=\,-\div \nabla F\,
=\,\div(\lambda\'\s),
\quad\text{by (iv),} \\
\vspace{1ex}
&=\,\<\nabla\lambda,\s\>\,+\,\lambda\div\s\,
=\,|\'\s\'|^2\,-\,n\'\lambda^2,
\quad\text{by (i),} \\
\vspace{1ex}
&=\,c^2-\lambda^2-n\'\lambda^2,
\quad\text{by (iii).}
\tag"\qed"
\endalign
$$
\enddemo
\proclaim{Theorem 5.2}
Let $\s$ be a conformal gradient field on $M=S^n$.  Then $\s$ is a
$(p,q)$-harmonic section of $TM$ if and only if $n\geqs3$ and:
$$
p=n+1,\quad
q=2-n,\quad 
\|\'\s\'\|_\infty=1/\sqrt{-q}.
$$
\endproclaim
\demo{Proof}
From Lemma 5.1 it follows that:
$$
|\'\nabla\s\'|^2\,=\,n\lambda^2,
\qquad
\nab{\nabla F}\s\,=\,\lambda^2\'\s,
\qquad
|\'\nabla F|^2\,=\,\lambda^2(c^2-\lambda^2)
$$
Therefore:
$$
\align
T_p(\s)\,&=\,(1+2F)\nabla^*\nabla\s\,
+\,2p\,\nab{\nabla F}\s \\
&=\,(1+c^2-\lambda^2)\s\,
+\,2p\'\lambda^2\s\,
=\,\bigl(1+c^2+(2p-1)\lambda^2\bigr)\s, 
\endalign
$$
$$
\align
\phi_{p,q}(\s)\,
&=\,p\,|\'\nabla\s\'|^2\,
-\,pq\,|\'\nabla F|^2\,
-\,q(1+2F)\Delta F \\
&=\,np\lambda^2\,
-\,pq\lambda^2(c^2-\lambda^2)\,
-\,q(1+c^2-\lambda^2)\bigl(c^2-(n+1)\lambda^2\bigr) \\
&=\,p(n+q)\lambda^2\,
-\,q(1+c^2-\lambda^2)\bigl(c^2-(n-p+1)\lambda^2\bigr)
\endalign
$$
Thus the harmonic section equations are polynomial in $\lambda$.  Since
$\lambda$ is a continuous function on $M$, which vanishes only on the great
hypersphere orthogonal to $a$, this polynomial is identically zero if and
only if the coeffiecients of the various powers of $\lambda$ vanish:
$$
\align
q\,&=\,-1/c^2,\quad\text{from the zeroth order terms,}
\tag5-1 \\
2p-1\,&=\,p(n+q)\,+\,qc^2\,+\,q(1+c^2)(n-p+1),
\quad\text{from $\lambda^2$,}
\tag5-2 \\
p\,&=\,n+1,\quad\text{from $\lambda^4$.}
\tag5-3
\endalign
$$
Substituting (5-1) and (5-3) into (5-2) yields:
$$
2(n+1)\,=\,q(n+1)\,+\,n(n+1),
$$
and hence $q=2-n$.
\qed
\enddemo
\subhead
Remark 5.3
\endsubhead
The $(p,q)$-harmonic sections $\s$ in Theorem 5.2 are clearly
$q$-Riemannian, but not $r$-Riemannian for any $r<q$.  Indeed a simple
calculation using Lemma 5.1 yields:
$$
|\'\nabla\s\'|^2\,+\,q\,|\'\nabla F|^2\,
=\,(n-1)\lambda^2+(n-2)\lambda^4,
$$
and this is strictly positive except when $\lambda=0$, which occurs
precisely on the great hypersphere orthogonal to the axial vector; in
particular $E^v(\s)\geqs0$.  However, since: 
$$
1-p\,=\,-n\,<\,2-n\,=\,q,
$$ 
it follows
that $\s$ is not $(1-p)$-Riemannian, and hence not $\nu(1-p)$-Riemannian
for any $\nu\geqs1$.  Therefore $\s$ is not amenable to Theorem 4.6.  In
fact $q<\varrho(p)$ if $n>3$, but $q=\varrho(p)$ when $M=S^3$.
\medskip
From Theorem 5.2, the only instance where standard conformal gradient fields
are $(p,q)$-harmonic sections of $TM$ is when $M=S^3$, in which case
$(p,q)=(4,-1)$.  The fact that in higher dimensions $(p,q)$-harmonic
sections are obtained by scaling the standard fields by a {\it constant\/}
factor (viz\. $1/\sqrt{n-2}\,$) is a manifestation of the non-linearity of
the $(p,q)$-harmonic section equations when $p\neq0$.  We now show that when
$n=2m+1$ the only {\it functional\/} multiples of the Hopf vector field
which are $(p,q)$-harmonic sections of $TM$ are those described in Example
4.2.
\proclaim{Theorem 5.4}
Suppose $\s=f\'\xi$, where $\xi$ is the Hopf vector field on $M=S^{2m+1}$
and $f\colon M\to\Bbb R$ is any smooth function.  Then $\s$ is a non-trivial
$(p,q)$-harmonic section of $TM$ if and only if $p>1$ and
$f=\pm1/\sqrt{p-1}\,$.
\endproclaim
\demo{Proof}
Initially, suppose $\xi$ is any unit vector field on any manifold $M$
with $\upchi(M)=0$.  Then:
$$
\align
2F\,&=\,|\'\s\'|^2\,
=\,|\'f\'\xi\'|^2\,=\,f^2,
\qquad
\nabla F\,
=\,\frac12\,\nabla(f^2)\,
=\,f\,\nabla f \\
\vspace{1ex}
\nab{\nabla F}\s\,
&=\,\nab{\nabla F}(f\'\xi)\,
=\,\bigl(\nabla F.f\bigr)\xi\,
+\,f\,\nab{\nabla F}\xi\,
=\,f\,|\'\nabla f\'|^2\'\xi\,
+\,f^2\,\nab{\nabla f}\xi \\
\vspace{1ex}
\nabla^*\nabla\s\,
&=\,\nabla^*\nabla(f\'\xi)\,
=\,f\,\nabla^*\nabla\xi\,
+\,(\Delta f)\'\xi\,-\,2\,\nab{\nabla f}\xi
\endalign
$$
It follows that:
$$
\align
T_p(\s)\,&=\,(1+2F)\nabla^*\nabla\s\,
+\,2p\,\nab{\nabla F}\s \\
\vspace{1ex}
&=\,(1+f^2)f\,\nabla^*\nabla\xi\,
+\,2\bigl((p-1)f^2-1\bigr)\nab{\nabla f}\xi\,
+\,\bigl(2pf\,|\'\nabla f\'|^2\,
+\,(1+f^2)\Delta f\bigr)\xi
\endalign
$$
We also have:
$$
\allowdisplaybreaks
\align
|\'\nabla\s\'|^2\,
&=\,\sum_i\bigl|\'\nab{E_i}(f\'\xi)\'\bigr|^2\,
=\,\sum_i\bigl|(E_i.f)\xi\,
+\,f\,\nab{E_i}\xi\'\bigr|^2 \\
&=\,\sum_i\left(|\'E_i.f\'|^2\,
+\,2f(E_i.f)\<\xi,\nab{E_i}\xi\>\,
+\,f^2\,|\'\nab{E_i}\xi\'|^2\right) \\
&=\,|\'\nabla f\'|^2\,
+\,f^2\,|\'\nabla\xi\'|^2,
\quad\text{since $|\'\xi\'|=1$,} \\
\vspace{1ex}
|\'\nabla F|^2\,&=\,f^2\,|\'\nabla f\'|^2, \\
\vspace{1ex}
\Delta F\,&=\,\frac12\,\Delta(f^2)\,
=\,f\,\Delta f\,-\,|\'\nabla f\'|^2, 
\endalign
$$
where to compute $\Delta(f^2)$ we have used formula (4-8).  It follows
that:
$$
\align
\phi_{p,q}(\s)\,
&=\,p\,|\'\nabla\s\'|^2\,
-\,q(1+2F)\'\Delta F\,
-\,pq\,|\'\nabla F|^2 \\
\vspace{1ex}
&=\,pf^2\,|\'\nabla\xi\'|^2\,
+\,\bigl(p+q+(1-p)qf^2\bigr)|\'\nabla f\'|^2\,
-\,q(1+f^2)f\'\Delta f
\endalign
$$
Therefore:
$$
\align
T_p(\s)&-\phi_{p,q}(\s)\'\s\,
=\,(1+f^2)f\,\nabla^*\nabla\xi\,
+\,2\bigl((p-1)f^2-1\bigr)\nab{\nabla f}\xi \\
\vspace{1ex}
&+\,\left(\vphantom{\nab{\nabla F}\s}
(1+f^2)(1+rf^2)\'\Delta f\,
+\,\bigl((p-1)rf^2+p-q\bigr)f\,|\'\nabla f\'|^2\,
-\,pf^3\,|\'\nabla\xi\'|^2\right)\xi
\endalign
$$
Now suppose that $\xi$ is a harmonic section of $S\E(1)$, and therefore
satisfies equation (1-3) with $k=1$.  Then:
$$
\align
T_p(\s)-\phi_{p,q}(\s)\'\s\,
&=\,2\bigl((p-1)f^2-1\bigr)\nab{\nabla f}\xi \\
\vspace{1ex}
&\quad+\,\left(\vphantom{\nab{\nabla F}\s}
(1+f^2)(1+qf^2)\'\Delta f\,
+\,\bigl((p-1)qf^2+p-q\bigr)f\,|\'\nabla f\'|^2\right. \\
&\qquad\quad\left.
+\;\bigl(1+(1-p)f^2\bigr)f\,|\'\nabla\xi\'|^2
\vphantom{\nab{\nabla F}\s}\right)\xi
\tag5-4
\endalign
$$
Since $\nabla\xi$ is orthogonal to $\xi$, because $|\'\xi\'|=1$, it follows
from (5-4) that a necessary condition for $\s$ to be $(p,q)$-harmonic is:
$$
\nab{\nabla f}\xi\,=\,0
\tag5-5
$$
Finally suppose that $\xi$ is the Hopf vector field on $M=S^{2m+1}$
(Example 4.2). Then:
$$
\nab X\xi\,
=\,\cases 
iX,&\text{if $X\perp\xi$,}\\
0,&\text{if $X\parallel\xi$,}
\endcases
\tag5-6
$$
which implies that (5-5) holds if and only if $\nabla f=\mu\'\xi$ for some
smooth $\mu\colon M\to\Bbb R$, or equivalently:
$$
df\,=\,\mu\,\xi^{\flat},
\tag5-7
$$ 
where $\xi^{\flat}$ is the $1$-form metrically dual to $\xi$.  It also
follows from (5-6) that $d\xi^{\flat}=\omega$, the restriction of the
K\"ahler form of $\Bbb C^{m+1}$, for:
$$
\align
2\'d\'\xi^{\flat}(X,Y)\,
&=\,\nab X\xi^{\flat}(Y)\,
-\,\nab Y\xi^{\flat}(X)\,
=\,\bigl\<\nab X\xi,Y\bigr\>\,
-\,\bigl\<\nab Y\xi,X\bigr\> \\
&=\,\<iX,Y\>\,-\,\<iY,X\>,
\quad\text{for \underbar{all} $X,Y$} \\
&=\,2\<iX,Y\>\,=\,2\omega(X,Y)
\endalign
$$
Exterior differentiation of (5-7) therefore yields the following
differential equation for
$\mu$:
$$
0\,=\,d\mu\wedge\xi^{\flat}\,+\,\mu\,\omega
\tag5-8
$$
In particular, if $A\perp\xi$ then $iA\perp\xi$ also and:
$$
\omega(A,iA)\,=\,\<iA,iA\>\,=\,|\'A\'|^2
$$
Now:
$$
2\'d\mu\wedge\xi^{\flat}(X,Y)\,
=\,d\mu(X)\'\<\xi,Y\>\,
-\,d\mu(Y)\'\<\xi,X\>,
$$
and so: 
$$
d\mu\wedge\xi^{\flat}(A,iA)\,=\,0
$$  
Therefore the only solution of (5-8) is $\mu=0$, and it follows that $f$ is
constant.  But then (5-4) reduces to:
$$
T_p(\s)-\phi_{p,q}(\s)\'\s\,
=\,\left(2m\bigl(1+(1-p)f^2\bigr)f
\vphantom{\nab{\nabla F}\s}\right)\xi,
$$
so if $\s$ is $(p,q)$-harmonic then either $f=0$ or $f^2=1/(p-1)$.
\qed
\enddemo

\enddocument
\end